\documentclass[reqno]{amsart}
\usepackage{amssymb, amscd}
\newtheorem{theorem}{Theorem}[section]
\newtheorem{lemma}[theorem]{Lemma}
\newtheorem{corollary}[theorem]{Corollary}
\theoremstyle{definition}
\newtheorem{remark}[theorem]{Remark}
\newtheorem{conj}{Conjecture}
\numberwithin{equation}{section}

\newcommand{\z}{\mathbb{Z}}
\newcommand{\maps}[3]{#1\colon #2\to #3} 
\newcommand{\fld}{\mathbb{F}} 
\newcommand{\Fld}{\mathbb{E}} 
\newcommand{\val}{\mathfrak{v}}
\newcommand{\cism}{\mathfrak{i}} 
\DeclareMathOperator{\Ker}{Ker}
\DeclareMathOperator{\rank}{rank}
\DeclareMathOperator{\rep}{Rep} 
\DeclareMathOperator{\Hom}{Hom}
\DeclareMathOperator{\sign}{sign} 
\DeclareMathOperator{\rad}{rad} 
\DeclareMathOperator{\ord}{ord} 

\begin{document}
\title{Quadratic quandles and their link invariants}

\author{R.A.~Litherland} 
\address{Department of Mathematics\\ 
         Louisiana State University \\ 
         Baton Rouge, LA 70803} 
\email{lither@math.lsu.edu} 
\date{July 11, 2002}
\subjclass{Primary: 57M25; Secondary 55N99}
\keywords{Knots, links, quandle cohomology, state-sum invariants}

\begin{abstract}
  Carter, Jelsovsky, Kamada, Langford and Saito have defined an
  invariant of classical links associated to each element of the
  second cohomology of a finite quandle. We study these invariants for
  Alexander quandles of the form 
  $\z[t,t^{-1}]/(p, t^2 +\kappa t + 1)$, where $p$ is a prime number
  and $t^2 + \kappa t + 1$ is irreducible modulo $p$. For each such
  quandle, there is an invariant with values in the group ring
  $\z[C_p]$ of a cyclic group of order $p$. We shall show that the
  values of this invariant all have the form $\Gamma_p^r p^{2s}$ for a
  fixed element $\Gamma_p$ of $\z[C_p]$ and integers $r \geq 0$ and 
  $s > 0$. We also describe some machine computations, which lead us
  to conjecture that the invariant is determined by the Alexander
  module of the link. This conjecture is verified for all torus and
  two-bridge knots.  
\end{abstract}

\maketitle

\section{Quadratic quandles}
In \cite{CJKLS}, Carter, Jelsovsky, Kamada, Langford and Saito
associate an invariant of classical links to each element of
the second cohomology of a finite quandle.  In this paper we study
these invariants for a class of quandles we refer to as quadratic
quandles, which will be described shortly. 
We assume familiarity with the definitions of quandle and quandle
cohomology as given in \cite{CJKLS}. Our notation generally follows
that paper, except that we write the operation in a quandle as
exponentiation (following Fenn and Rourke \cite{FR}).
We denote the ring $\z[t,t^{-1}]$ of integer Laurent polynomials by
$\Lambda$, and the ring $\z/(n)[t,t^{-1}]$ by $\Lambda_n$.  For any
ring $R$, we denote its group of units by $R^\times$.

Recall that an Alexander quandle is a $\Lambda$-module with the
quandle operation $a^b=ta+(1-t)b$.  We let $p$ be a prime number, and
denote the finite fields of orders $p$ and $p^2$ by $\fld_p$ and
$\Fld_p$ respectively. For $x \in \Fld_p$, we denote the conjugate of
$x$ over $\fld_p$ by $\bar x$ and its norm $x\bar x \in \fld_p$ by
$N(x)$.  We fix a generator $\theta$ of $\Fld_p$ over $\fld_p$, and
give $\Fld_p$ the structure of a $\Lambda$-module by setting $tx =
\theta x$ ($x \in \Fld_p$). Thus $\Fld_p \cong
\Fld_p[t,t^{-1}]/(t-\theta)$ as a $\Lambda$-module, so Theorem 2.2 of
Mochizuki \cite{M} gives the dimension of the quandle cohomology
$H^2_Q(\Fld_p;\Fld_p)$ over $\Fld_p$, which is the same as that of
$H^2_Q(\Fld_p;\fld_p)$ over $\fld_p$; the dimension is 1 if
$\theta^{p+1} = 1$, and 0 otherwise. Since the non-trivial element of 
the Galois group of $\Fld_p$ over $\fld_p$ is the Frobenius
automorphism, $\theta^{p+1} = \theta\bar\theta= N(\theta)$. Hence we
further require $\theta$ to have norm 1, and when this is so we call
$\Fld_p$ with the structure of an Alexander quandle determined by
$\theta$ a quadratic quandle. We denote the minimal polynomial of
$\theta$ over $\fld_p$ by $h$; it is an irreducible polynomial of the
form $t^2+\kappa t+1$. An equivalent description is that a quadratic
quandle is an Alexander quandle $\Lambda_p/(h)$ where $h$ is such a
polynomial. 

For the quadratic quandle $\Fld_p$ we have $H^2_Q(\Fld_p;\fld_p) \cong
\fld_p$, and choosing a non-zero element gives a link invariant with 
values in the group ring $\z[C_p]$ of a cyclic group of order $p$. To
describe an explicit cocycle, we fix a non-zero element  
$\cism\in \Fld_p$ with $\bar \cism = -\cism$. Define 
$\maps \phi{\Fld_p^2}{\fld_p}$ by $\phi(x, y) = 
\cism(x\bar y - \bar x y)$. Then $\phi$ is an alternating
$\fld_p$-bilinear form, and satisfies 
\begin{equation}\label{eq:Nphi}
  \phi(zx, zy) = N(z)\phi(x, y)\quad\text{for $x, y, z \in \Fld_p$.}
\end{equation}
Recall that the quandle chain group $C_2^Q(\Fld_p)$ is obtained from
the free abelian group on $\Fld_p^2$ by quotienting out the subgroup
generated by pairs $(x,x)$. Since $\phi(x,x)=0$, $\phi$ determines a
cochain in $C^2_Q(\Fld_p;\fld_p)$, which we also denote by $\phi$.

\begin{lemma}\label{lem:cocycle} 
  For the quadratic quandle $\Fld_p$, the function $\phi$ is a cocycle
  representing a non-zero element of $H^2_Q(\Fld_p;\fld_p)$. 
\end{lemma}

\begin{proof}
  Let $x, y, z \in \Fld_p$. From \eqref{eq:Nphi} we have 
  $\phi(\theta x, \theta y) = \phi(x, y)$ and so  
  \[
    \phi\bigl(\theta x, (1-\theta)y\bigr) 
    = \phi(\theta x, y) -\phi(x, y) = \phi\bigl((\theta-1)x, y\bigr).
  \]
  Hence 
  \begin{align*}
    \phi(x^z, y^z) 
    &= \phi\bigl(\theta x + (1-\theta)z, \theta y + (1-\theta)z\bigr)\\
    &=\phi(\theta x, \theta y) + \phi\bigl(\theta x, (1-\theta)z\bigr) +
    \phi\bigl((1-\theta)z, \theta y\bigr)\\ 
    &=\phi(x, y) + \phi\bigl(\theta(x-y), (1-\theta)z\bigr)\\ 
    &=\phi(x,y) + \phi\bigl((\theta-1)(x-y), z\bigr)\\ 
    &=\phi(x, y) + \phi\bigl(\theta x + (1 - \theta)y, z\bigr)-\phi(x,z)\\ 
    &=\phi(x, y) + \phi(x^y,z)-\phi(x,z),
  \end{align*}
  so
  \[
    \phi\partial(x, y, z) 
    = \phi(x,z) - \phi(x^y,z) - \phi(x, y) +\phi(x^z,y^z)
    = 0. 
  \]
  Thus $\phi$ is a cocycle.  To see that it is not a coboundary, we
  observe that $c = (1, -\theta) - (1, 0) - (\theta, 0) \in
  C_2^Q(\Fld_p)$ is a cycle and $\phi(c) = \cism(\theta - \bar\theta)
  \neq 0$. 
\end{proof}

\section{Their link invariants}
Let $K$ be an oriented classical link, and $Q(K)$ its fundamental
quandle. The cocycle invariant $\Phi_\phi(K)$ defined in \cite{CJKLS}
involves a sum over all quandle homomorphisms from $Q(K)$ to
$\Fld_p$. The number of such homomorphisms is determined in Theorem 1
of Inoue \cite{Inou}. The proof of that result may be viewed in the
following way. (See Example 2 in \S6.4 of \cite{FR}.) We first recall
some facts about the (one-variable) Alexander module of $K$. Let $X$
be the exterior of $K$, $\tilde X$ its infinite cyclic cover, and
$\tilde X_0$ a fiber of $\tilde X$. Then the Alexander module is
$H_1(\tilde X, \tilde X_0)$, regarded as a $\Lambda$-module with $t$
acting by a generator of the group of covering transformations. We
shall use the notation $A(K)$ for $H_1(\tilde X, \tilde X_0)$.  We
also set $\hat A(K)= H_1(\tilde X)$, and then there is a short exact
sequence
\[
  0 \to \hat A(K) \to A(K) \to \Lambda \to 0.
\]
(If $K$ is a knot, $\hat A(K)$ is the torsion submodule of $A(K)$.)
Similarly, we have $\Lambda_p$-modules $A_p(K) = 
H_1(\tilde X, \tilde X_0; \fld_p)$ and $\hat A_p(K)= 
H_1(\tilde X;\fld_p)$.

Now take a diagram $D$ of $K$. At each crossing $c$ let $\omega(c)$ be
the overcrossing arc, $\lambda(c)$ the undercrossing arc to the left
of $\omega(c)$, and $\rho(c)$ the one to the right.  Then $A(K)$ has a
$\Lambda$-module presentation with the arcs of $D$ as generators and a
relation $\lambda(c) = t\rho(c) + (1-t)\omega(c)$ for each
crossing. Also, the fundamental quandle of $K$ has a presentation with
the arcs of $D$ as generators and a relation $\lambda(c) =
\rho(c)^{\omega(c)}$ for each crossing.  Thus, when $A(K)$ is regarded
as a quandle, there is a quandle homomorphism $Q(K) \to A(K)$ which is
universal among homomorphisms from $Q(K)$ to Alexander quandles. Hence
there is a one-to-one correspondence between quandle homomorphisms
$Q(K) \to \Fld_p$ and $\Lambda$-module homomorphisms $A(K) \to
\Fld_p$, which in turn correspond to $\Lambda_p$-module homomorphisms
$A_p(K) \to \Fld_p$. We set $\rep(K, \Fld_p) =
\Hom_{\Lambda_p}(A_p(K), \Fld_p)$; this is a vector space over
$\Fld_p$. Now $\Lambda_p$ is a PID, so we may consider the primary
decomposition of $A_p(K)$. Let $\nu_h(K)$ be the number of $h$-torsion
factors in this decomposition, and let $\nu_0(K)$ be the rank of
$A_p(K)$. (If $K$ is a knot, $\nu_0(K) = 1$.) Then the dimension of
$\rep(K, \Fld_p)$ over $\Fld_p$ is $\nu_h(K) + \nu_0(K)$. We shall
sometimes regard an element of $\rep(K,\Fld_p)$ as a coloring of the
arcs of $D$ with elements of $\Fld_p$ satisying the appropriate
relation at each crossing. We also denote the $h$-torsion part of
$A_p(K)$ by $A_p(K)_h$.

The cocycle invariant may be defined as follows. (Our terminology and
notation differs slightly from that of \cite{CJKLS}, but the
definitions are equivalent.) Let $C_p$ be a (multiplicative) cyclic
group of order $p$, with generator $u$. Let $\sign(c) = \pm 1$ be the
sign of the crossing $c$. For each $f\in \rep(K, \Fld_p)$, define the
Boltzmann weight of $f$ to be
\[
  B(f) = \sum_c \sign(c)\phi\bigl(f\rho(c), f\omega(c)\bigr) \in\fld_p, 
\]
where the sum is taken over all crossings of the diagram $D$. Then
\[
  \Phi_\phi(K) = \sum_{f \in \rep(K, \Fld_p)}u^{B(f)} \in \z[C_p].
\]
The invariant depends only on the cohomology class of $\phi$; for
quadratic quandles, it follows from Theorem \ref{thm:main} below that
all non-zero classes give the same invariant.  Now define
$\maps\eta{\rep(K, \Fld_p) \times \rep(K, \Fld_p)}{\Fld_p}$ by
\[
  \eta(f,g) = \sum_c \sign(c)\cism\bigl(f\rho(c)\overline{g\omega(c)}
  -f\omega(c)\overline{g\rho(c)}\bigr).
\]
Then $\eta$ is a Hermitian form with respect to the involution 
$x \mapsto \bar x$ of $\Fld_p$, and $B(f) = \eta(f, f)$ for $f \in
\rep(K,\Fld_p)$.
\begin{theorem} \label{thm:main} 
  Let $\Gamma_p = 1+(p+1)\sum_{i=1}^{p-1}u^i \in \z[C_p]$. For the
  cocycle of Lemma \ref{lem:cocycle}, we have $\Phi_\phi(K) = 
  \Gamma_p^r p^{2s}$, where $r$ is the rank of the Hermitian form $\eta$
  and $s = \nu_h(K) + \nu_0(K) - r$ is its nullity. Further, $s > 0$.
\end{theorem}
\begin{proof} 
  Let $n = \nu_h(K) + \nu_0(K)$, the dimension of $\rep(K,\Fld_p)$. We
  may choose a basis $\{f_1,\ldots,f_n\}$ of $\rep(K,\Fld_p)$ that is
  orthogonal with respect to $\eta$. Let $a_i = \eta(f_i,f_i) \in
  \fld_p$. Then, for $\vec x = (x_1,\ldots,x_n) \in \Fld_p^n$, we have
  $B(\sum_{i=1}^n x_i f_i) = \sum_{i=1}^n a_iN(x_i)$, and so
  \[
    \Phi_\phi(K) = \sum_{\vec x \in \Fld_p^n}u^{\sum_{i=1}^n a_iN(x_i)}
    = \prod_{i=1}^n\sum_{x\in \Fld_p}u^{a_iN(x)}.
  \]
  Since the norm restricts to a homomorphism of $\Fld_p^\times$ onto
  $\fld_p^\times$, $\sum_{x\in \Fld_p}u^{a_iN(x)}$ is equal to
  $\Gamma_p$ if $a_i\neq 0$, and $p^2$ if $a_i = 0$. It only remains
  to prove that $s > 0$; that is, that the form $\eta$ is
  singular. This we do in Lemma \ref{lem:sing} below.  
\end{proof}
\begin{lemma} \label{lem:rel} 
  In the Alexander module of a link $K$ with diagram $D$, there is a
  relation
  \[
    (1-t)\sum_c \sign(c)\bigl(\rho(c) - \omega(c)\bigr) = 0.
  \]
\end{lemma}
\begin{proof} 
  For each crossing $c$, let $\alpha(c)$ be the undercrossing arc
  oriented into $c$, and $\beta(c)$ the one oriented out. Thus, if
  $\sign(c) = 1$ then $\alpha(c) = \rho(c)$ and $\beta(c) =
  \lambda(c)$, while if $\sign(c) = -1$ then $\alpha(c) = \lambda(c)$
  and $\beta(c) = \rho(c)$. We have  
  \begin{align*}
    (1-t)\sum_c\sign(c)\bigl(\rho(c) - \omega(c)\bigr) &= 
    \sum_c\sign(c)\bigl(\rho(c) - \lambda(c)\bigr)\\
    &= \sum_c \bigl(\alpha(c) - \beta(c)\bigr)\\ 
    &= 0,
  \end{align*}
  because each arc of $D$ is $\alpha(c)$ for a single crossing $c$,
  and $\beta(c')$ for another crossing $c'$. 
\end{proof}
\begin{remark} 
  If $K$ is a knot, then $A(K)$ has no $(1-t)$-torsion, so we have the
  relation $\sum_c \sign(c)\bigl(\rho(c) - \omega(c)\bigr) = 0$. 
\end{remark}
Let $\rep_0(K,\Fld_p) \leq \rep(K, \Fld_p)$ be the annihilator of
$\hat A_p(K)$, a 1-dimensional subspace of $\rep(K, \Fld_p)$. Because
all arcs of $D$, considered as elements of $A_p(K)$, map to the same
element of $A_p(K)/\hat A_p(K)$, each element of $\rep_0(K,\Fld_p)$
takes the same value on each arc. Thus elements of $\rep_0(K, \Fld_p)$
correspond to constant colorings of $D$, and have Boltzmann weight 0;
however, more is true. We denote the radical of $\eta$ by
$\rad\eta$. Thus $g\in \rep(K, \Fld_p)$ is in $\rad \eta$ iff 
$\eta(f,g) = 0$ for all $f \in \rep(K, \Fld_p)$. 
\begin{lemma}\label{lem:sing} 
  $\rep_0(K, \Fld_p)$ is contained in  $\rad\eta$; in particular,
  $\eta$ is singular. 
\end{lemma}
\begin{proof} 
  Let $f \in \rep(K, \Fld_p)$ and $g \in \rep_0(K, \Fld_p)$, and let
  $x \in \Fld_p$ be the common value of $g$ on the arcs of $D$. Then  
  \[
    \eta(f, g) =\sum_c \sign(c)\cism\bigl(f\rho(c)-f\omega(c)\bigr)\bar x
    =\cism\bar xf(a),
  \]
  where $a = \sum_c \sign(c)\bigl(\rho(c) - \omega(c)\bigr) \in A_p(K)$.
  By Lemma \ref{lem:rel}, 
  \[
    (1 -\theta)f(a) =f\bigl((1-t)a\bigr) = 0.
  \]
  Since $1-\theta$ is a non-zero element of the field $\Fld_p$, $f(a)
  = 0$ and we are done.  
\end{proof}
Note that if $\nu_0(K)=1$, this lemma implies that
$\rank\eta\leq\nu_h(K)$. This is the case for all knots, and for all
links considered in the rest of this paper.

\section{Computations}
We have carried out some machine computations of the invariants of
knots associated to quadratic quandles. The program (written in C)
proceeds as follows. The knots are taken from the tables distributed
with Knotscape \cite{HT}, which describe minimal-crossing diagrams in
Dowker-Thistlethwaite code. From this, a presentation matrix for
$A_p(K)$ is read off. This matrix is diagonalized, and the matrix
expressing the original (arc) generators in terms of the new
generators is computed. We find $A_p(K)_h$ by finding the highest
power of $h$ dividing each non-zero entry; in particular, $\nu_h(K)$
is just the number of non-zero entries divisible by $h$.  The form
$\eta'(f,g) = \eta(f,g) + \eta(g,f)$ is a symmetric $\fld_p$-bilinear
form on $\rep(K, \Fld_p)$ with $\rank\eta'= 2\rank \eta$.  Using the
change-of-basis matrix, a matrix for $\eta'$ can be found, and its
rank determined, which gives $\rank\eta$. (It may seem inefficient to
use $\eta'$ instead of $\eta$ since this doubles the size of the
matrix, but it is easier to work over $\fld_p$ rather than $\Fld_p$.)
Now $\Phi_\phi(K)$ is determined by Theorem \ref{thm:main}. Some
computations of quandle invariants are given by Carter, Jelsovsky,
Kamada and Saito \cite{CJKS}, and there is some overlap with ours. For
$p=2$ or 3 there is just one choice for $h$, giving two quadratic
quandles $\Lambda_2/(t^2+t+1)$ and $\Lambda_3/(t^2+1)$.  In
\cite{CJKS}, the values of the invariants for these quandles are
computed for knots up to 9 crossings, and we have repeated these
computations. For $p=2$, the results agree, but for $p=3$ there are
some discrepancies, at least one of which is due to a typographical
error in \cite{CJKS}. To obtain our results from theirs: 
\begin{enumerate}
  \item in the list of knots with invariant
    $9(1+4t+4t^2),$\footnote{Their $t$ is our $u$.} replace $9_{47}$
    by $9_{48}$;  
  \item from the list of knots with invariant 81, delete $8_{24}$
    [sic] and $9_{46}$, and add $8_6$, $9_3$ and $9_{47}$. 
\end{enumerate}
We have not extended the tabulation of invariants to higher crossing
numbers; instead, we have counted the knots with a given combination
of $A_p(K)_h$ and $\rank \eta$. For example, with $p = 2$, of the
313230 knots of at most 15 crossings, 445 have $\nu_h(K)$ equal to
3. These may be further classified as in Table \ref{tab:nh3}, where an
entry $(a,b,c)$ in the first column means that $A_2(K)_h\cong
\Lambda_2/(h^a)\oplus\Lambda_2/(h^b)\oplus\Lambda_2/(h^c)$.  
\begin{table}
  \begin{center}
    \begin{tabular}{|ccc|} 
      \hline 
      $A_2(K)_h$\rule{0pt}{11pt} & $\rank\eta$ & knots \\[2pt] 
      \hline
      $(1,1,1)$\rule{0pt}{11pt}& 3 & 217\\ 
      $(1,1,2)$ & 2 & 114\\ 
      $(1,1,3)$ & 2 & 110\\ 
      $(1,1,4)$ & 2 & \phantom{00}2\\ 
      $(1,2,2)$ & 1 & \phantom{00}1\\ 
      $(1,2,3)$ & 1 & \phantom{00}1\\[2pt] 
      \hline
    \end{tabular}
    \vspace{15pt} 
    \caption{}
    \label{tab:nh3}
  \end{center}
\end{table}
It will be noted that in these cases $\rank\eta$ is equal to the
number of ones in the triple $(a,b,c)$. In general, let $\nu_h'(K)$ be
the number of factors in the primary decomposition of $A_p(K)$ of the
form $\Lambda_p/(h)$. 
\begin{conj}\label{conj:weak}
  For any link $K$, $\rank\eta=\nu'_h(K)$.
\end{conj}
If this is true, the invariants associated to quadratic quandles
contain no more information than the Alexander module, which would be
somewhat disappointing. We have verified it for $p = 2$, 3, 5, 7 and
11, for each possible choice of the coefficient $\kappa$ in $h$
(giving twelve invariants), and for all knots up to (and including) 16
crossings. In fact, a little more seems to be true. The number
$\nu_h'(K)$ has a more conceptual description.  Let $A_p'(K)$ be the
submodule of all $a \in A_p(K)$ for which $ha = 0$, and let
$\rep'(K,\Fld_p) \leq \rep(K, \Fld_p)$ be its annihilator. 
\begin{lemma}\label{lem:repp}
  The codimension of $\rep'(K, \Fld_p)$ in $\rep(K,\Fld_p)$ is
  $\nu_h'(K)$. 
\end{lemma}
\begin{proof} 
  Suppose that $A_p(K)_h$ is isomorphic to
  $\bigoplus_{i=1}^{\nu_h(K)}\Lambda_p/(h^{e_i})$, and let $a_i$ be a
  generator of the $i^{\rm th}$ factor. Then $A_p'(K)$ is generated by
  the elements $h^{e_i-1}a_i$. If $e_i > 1$, then $h^{e_i-1}a_i$ is
  annihilated by any element of $\rep(K, \Fld_p)$, so $f \in
  \rep(K,\Fld_p)$ belongs to $\rep'(K,\Fld_p)$ iff $f(a_i) = 0$ for
  those values of $i$ for which $e_i = 1$. 
\end{proof}
\begin{conj}\label{conj:strong}
  For any link $K$, $\rad\eta = \rep'(K, \Fld_p)$.
\end{conj}
This has been verified for the same range of invariants and knots as
Conjecture \ref{conj:weak}. In the rest of the paper, we shall see
that it is true for all torus and two-bridge knots. In most cases, we
need only check Conjecture \ref{conj:weak}, because of the following
simple observation. 
\begin{lemma}\label{lem:c1impc2}
  Suppose that $\rank\eta = \nu_h'(K)$, and that either
  $\nu_h'(K)=\nu_h(K)$ and $\nu_0(K)=1$, or $\nu_h'(K)=0$. Then
  $\rad\eta = \rep'(K, \Fld_p)$.
\end {lemma}
\begin{proof}
  In the first case we have $\rad\eta=\rep_0(K, \Fld_p)
  =\rep'(K,\Fld_p)$, and in the second, $\rad\eta=
  \rep(K, \Fld_p)=\rep'(K,\Fld_p)$. 
\end{proof}

\section{Torus links}
If $R$ is a UFD and $\mathfrak{p}$ is a prime of $R$, we denote the
discrete $\mathfrak{p}$-adic valuation on $R$ by
$\val_{\mathfrak{p}}$. 
\begin{lemma}\label{lem:val}
  Let $\mathbb{K}$ be a field of characteristic $p$, let $a$ be a
  non-zero element of $\mathbb{K}$, let $f$ be an irreducible element
  of $\mathbb{K}[t,t^{-1}]$, and let $n$ be a positive integer. If
  $f\mid t^n-a^n$ then $\val_f(t^n-a^n) = p^{\val_p(n)}.$
\end{lemma}
\begin{proof}
  Let $e = \val_p(n)$, so that $n = p^e m$ where $p \nmid m$. Since
  $t^m-a^m$ is coprime to its derivative, it is square-free. Since
  $t^n-a^n = (t^m-a^m)^{p^e}$, the result follows.
\end{proof}
One case of this lemma we need is when $\mathbb{K} = \fld_p$, $a=1$,
and $f=h$. We denote the order of $\theta$ in $\Fld_p^\times$ by
$q$. Then $h \mid t^n - 1$ iff $q \mid n$.   

Let $T_{m,n}$ be the $(m, n)$ torus link.  If $\bar K$ is the mirror
image of a link $K$, $\Phi_\phi(\bar K) = \Phi_{\bar \phi}(K)$, where
$\bar \phi$ is a cocycle for the dual quandle to $\Fld_p$. Since
quadratic quandles are self-dual (because $h$ is a symmetric
polynomial), we may as well assume that $m$ and $n$ are both
positive. Let $c$ be the number of components of $T_{m,n}$, which is
the highest common factor of $m$ and $n$. We first determine the
Alexander module of $T_{m,n}$.
\begin{lemma} 
  If $c=1$, then $\hat A(T_{m,n}) \cong
  \Lambda/\bigl(\frac{(t^{mn}-1)(t-1)}{(t^m -1)(t^n-1)}\bigr)$. If 
  $c \geq 2$, $\hat A(T_{m,n})$ is isomorphic to the direct sum of
  $c-2$ copies of $\Lambda/(t^{mn/c}-1)$ and a module with
  presentation matrix 
  \[
    \left[\begin{matrix} 
      \frac{t^{mn/c}-1}{t^m - 1} & \frac{t^{mn/c}-1}{t^n - 1}\\ 
      0                          & (t-1)\frac{t^{mn/c}-1}{t^n - 1}
    \end{matrix}\right].
  \] 
  (Here the columns correspond to generators and the rows to relations.)
\end{lemma}
\begin{proof} 
  The exterior $X$ of $T_{m,n}$ is the union of two solid tori $V_1$
  and $V_2$, whose cores have linking numbers $m$ and $n$ respectively
  with $T_{m,n}$, and whose intersection consists of $c$ annuli. Let
  $\tilde V_1$ and $\tilde V_2$ be their inverse images in the
  infinite cyclic cover $\tilde X$. Then $H_1(\tilde V_1) = 
  H_1(\tilde V_2) = 0$, so $\hat A(T_{m,n})$ is isomorphic to the
  kernel of $H_0(\tilde V_1 \cap \tilde V_2) \to H_0(\tilde V_1)
  \oplus H_0(\tilde V_2)$. Further, $H_0(\tilde V_1) \cong \Lambda
  /(t^m-1)$, $H_0(\tilde V_2) \cong \Lambda /(t^n - 1)$, and
  $H_0(\tilde V_1 \cap \tilde V_2) \cong  
  \bigl(\Lambda/(t^{mn/c} -1)\bigr)^c$, and we may compute the above
  map as follows. Take points $x_1,\ldots,x_c$, one from each
  component of $V_1\cap V_2$, and for $2 \leq i\leq c$ take paths
  $\sigma_i$ from $x_1$ to $x_i$ in $V_1$ and $\tau_i$ from $x_i$ to
  $x_1$ in $V_2$. This may be done in such a way that $\sigma_i\tau_i$
  has linking number $i-1$ with $T_{m,n}$. Let $\tilde x_1$ be in the
  inverse image of $x_1$, and for $2 \leq x \leq c$ let
  $\tilde\sigma_i$ be the lift of $\sigma_i$ starting at $\tilde x_1$,
  let $\tilde x_i$ be its terminal point, and let $\tilde \tau_i$ be
  the lift of $\tau_i$ starting at $\tilde x_i$. Note that the
  terminal point of $\tilde \tau_i$ is $t^{i-1}\tilde x_1$. We take
  $\tilde x_1,\ldots\tilde x_c$ as generators for 
  $H_0(\tilde V_1 \cap \tilde V_2)$, and $\tilde x_1$ as a generator
  of each of $H_0(\tilde V_1)$ and $H_0(\tilde V_2)$.  Each $\tilde
  x_i$ is homologous to $\tilde x_1$ in $\tilde V_1$ (by 
  $\tilde \sigma_i$), and to $t^{i-1}\tilde x_1$ in $\tilde V_2$ (by
  $\tilde \tau_i$).  Thus we have a commutative diagram
  \[
    \begin{CD} 
      \Lambda^c                 @>\gamma>> \Lambda \oplus \Lambda \\
      @VV\alpha V                          @VV\beta V \\ 
      H_0(\tilde V_1 \cap \tilde V_2) @>>> H_0(\tilde V_1) 
      \oplus H_0(\tilde V_2)
    \end{CD}
   \] 
  in which $\alpha$ and $\beta$ are determined by our choice of
  generators and $\gamma(f_1,\ldots,f_c) = \bigl(\sum_{i=1}^c f_i,
  \sum_{i=1}^c t^{i-1}f_i\bigr)$. The case $c = 1$ follows easily (and
  is anyway well-known). Suppose that $c \geq 2$.  Then
  $(f_1,\ldots,f_c) \in\Ker(\beta\gamma)$ iff there exist $g_1,g_2\in
  \Lambda$ such that 
  \[
    \sum_{i=1}^c f_i = (t^m-1)g_1 \quad\text{and}\quad 
    \sum_{i=1}^c t^{i-1}f_i = (t^n-1)g_2.
  \]
  These equations imply that 
  \begin{align*}
    f_1 &= t\frac{t^m-1}{t - 1}g_1
          - \frac{t^n-1}{t-1}g_2 +\sum_{i=3}^c\frac{t^{i-1}-t}{t-1}f_i\\
    \text{and}\quad f_2 &= -\frac{t^m-1}{t - 1}g_1 + \frac{t^n-1}{t-1}g_2
          -\sum_{i=3}^c\frac{t^{i-1}-1}{t-1}f_i.
    \end{align*}
  Hence there is a surjection $\maps
  \delta{\Lambda^c}{\Ker(\beta\gamma)}$ defined by
  $\delta(g_1,g_2,f_3,\ldots,f_c) = (f_1,f_2,f_3,\ldots,f_c)$, with
  $f_1$ and $f_2$ given by the above formulas. Further, we have
  $(g_1,g_2,f_3,\ldots,f_c) \in \Ker(\alpha\delta)$ iff $f_i \equiv 0
  \pmod{t^{mn/c}-1}$ for $3 \leq i \leq c$ and there exist $h_1, h_2
  \in \Lambda$ such that 
  \begin{align*}
    t\frac{t^m-1}{t - 1}g_1 - \frac{t^n-1}{t-1}g_2 &= h_1(t^{mn/c}-1)\\  
    \text{and}\quad -\frac{t^m-1}{t - 1}g_1 + \frac{t^n-1}{t-1}g_2 &=
    h_2(t^{mn/c}-1). 
  \end{align*}
  But then
  \[
    g_1 = (h_1+h_2)\frac{t^{mn/c}-1}{t^m-1} \quad\text{and}
    \quad g_2 =(h_1+th_2)\frac{t^{mn/c}-1}{t^n-1}.
  \]
  The rows of the matrix in the statement of the lemma are given by
  the cases $h_1 = 1$, $h_2 = 0$ and $h_1 = -1, h_2 = 1$, so the
  result follows.
\end{proof}
Note that $\hat A_p(T_{m,n})$ is torsion, so the dimension of
$\rep(T_{m,n}, \Fld_p)$ is $\nu_h(T_{m,n}) + 1$, and $\rank\eta\leq
\nu_h(T_{m,n})$. 
\begin{lemma} 
  If $c=1$, then $\hat A_p(T_{m,n}) \cong
  \Lambda_p/\bigl(\frac{(t^{mn}-1)(t-1)}{(t^m -1)(t^n-1)}\bigr)$. If
  $c \geq 2$, $\hat A_p(T_{m,n})$ has the same $h$-torsion as
  $\Lambda_p/\bigl(\frac{t^{mn/c}-1}{t^m - 1}\bigr) \oplus
  \Lambda_p/\bigl(\frac{t^{mn/c}-1}{t^n - 1}\bigr) \oplus
  \bigl(\Lambda_p/(t^{mn/c}-1)\bigr)^{c-2}$.  
\end{lemma}
\begin{proof} 
  The case $c=1$ is obvious. If $c \geq  2$, we may assume (by
  symmetry in $m$ and $n$) that $p\nmid \frac{n}{c}$. Then $t-1$ and
  $\frac{t^{mn/c} - 1}{t^m - 1}$ are coprime in $\Lambda_p$ (by Lemma
  \ref{lem:val}), and it follows that $\hat A_p(T_{m,n})$ is isomorphic
  to the direct sum of $\Lambda_p/\bigl(\frac{t^{mn/c}-1}
  {t^m - 1}\bigr)$, $\Lambda_p/\bigl((t-1)\frac{t^{mn/c}-1}
  {t^n - 1}\bigr)$ and $c-2$ copies of $\Lambda_p/(t^{mn/c}-1)$.
  Since $\Lambda_p/\bigl((t-1)\frac{t^{mn/c}-1}{t^n - 1}\bigr)$ has
  the same $h$-torsion as $\Lambda_p/\bigl(\frac{t^{mn/c}-1}
  {t^n - 1}\bigr)$, we are done. 
\end{proof}
\begin{lemma}\label{lem:cases}\hfill
  \begin{itemize}
    \item[(0)] If $q \nmid \frac{mn}{c}$ then $\nu_h(T_{m,n}) =
      \nu'_h(T_{m,n}) = 0$. 
    \item[(1)] If $q \mid \frac{mn}{c}$, $p \mid c$ and $p \nmid
      \frac{m}{c}$ then  
      \begin{align*}
        \nu_h(T_{m,n}) &= \begin{cases}
          c,& \text{if $q \nmid n$ and either $q \nmid m$ or 
              $p \mid \frac{n}{c}$;}\\   
          c-2,& \text{if $q \mid m$, $q \mid n$ and $p \nmid \frac{n}{c}$;}\\ 
          c-1, &\text{otherwise}
        \end{cases}\\
        \text{and}\quad \nu'_h(T_{m,n}) &= 0.
      \end{align*}
    \item[(2)] If $q \mid \frac{mn}{c}$, $p \mid m$ and $p \nmid n$ then
      \begin{align*} 
        \nu_h(T_{m,n}) &= \begin{cases}
          c,& \text{if $q \nmid m$;}\\ 
          c-1, &\text{if $q \mid m$}
        \end{cases}\\
        \text{and} \quad \nu'_h(T_{m,n}) &=\begin{cases}
          1,& \text{if $p=2$, $3 \mid n$ and ${m}\equiv 2 \pmod 4$;}\\ 
          0,& \text{otherwise.}
        \end{cases}
      \end{align*}
    \item[(3)] If $q \mid \frac{mn}{c}$ and $p\nmid mn$ then
      \[
        \nu_h(T_{m,n}) = \nu'_h(T_{m,n}) = \begin{cases}
          c,& \text{if $q \nmid m$ and $q \nmid n$;}\\ 
          c-2,& \text{if $q \mid m$ and $q \mid n$;}\\ 
          c-1, & \text{otherwise.}
        \end{cases}
      \]
  \end{itemize}
\end{lemma}
Note that by symmetry in $m$ and $n$ (and the fact that $p$ cannot
divide both $\frac{m}{c}$ and $\frac{n}{c}$), we may always assume
that we are in one of the four cases of this lemma. The condition
$3\mid n$ appears in case 2 because $q=3$ when $p= 2$. 
\begin{proof} 
  By the previous lemma, $A_p(T_{m,n})_h\cong \bigoplus_{i=1}^c
  \Lambda/(h^{e_i})$, where if $c=1$ then $e_1 = 
  \val_h(t^{mn}-1) - \val_h(t^m-1)-\val_h(t^n-1)$, while if $c \geq 2$
  then $e_1 =\val_h(t^{mn/c}-1) - \val_h(t^m-1)$, $e_2 =
  \val_h(t^{mn/c}-1) - \val_h(t^n-1)$ and $e_i = \val_h(t^{mn/c}-1)$
  for $i \geq 3$. The number of non-zero $e_i$ is $\nu_h(T_{m,n})$,
  and the number equal to 1 is $\nu'_h(T_{m,n})$.  The rest of the
  proof is a simple, if tedious, case-by-case check using Lemma
  \ref{lem:val}.
\end{proof}
If $a$ and $b$ are coprime integers, let $\ord_a(b)$ denote the order
of $b$ in $\z/(a)^\times$.
\begin{theorem} \label{thm:torus} 
  For a torus link, $\rank\eta =\nu'_h(T_{m,n})$, except perhaps if 
  $q \mid \frac{mn}{c}$, $p \nmid mn$ and $4 \mid \ord_c(p)$. 
\end{theorem}
For a torus knot, the theorem applies for any $p$. More generally,
this is true if every odd prime factor of $c$ is congruent to $3 \pmod 
4$ and $16\nmid c$, since then $\z/(c)^\times$ has no
elements of order 4. For $p=2$, the smallest link not covered by the
theorem is $T_{5,15}$. In \cite[Remark 4.1]{CJKS}, the invariant of
this link is given as $544+480u = 4(1+3u)^4$. Thus $\rank\eta = 4$,
which is in fact equal to $\nu'_h(T_{5,15})$.  
\begin{corollary} 
  For a torus link, $\rad\eta =\rep'(T_{m,n},\Fld_p)$, except perhaps
  if $q \mid \frac{mn}{c}$, $p \nmid mn$ and $4 \mid \ord_c(p)$. 
\end{corollary}
\begin{proof} 
  This follows from Theorem \ref{thm:torus} and Lemma
  \ref{lem:c1impc2} except when $p=2$, $3 \mid n$, $m\equiv 2\pmod 4$,
  and $2 \nmid n$; this case will be dealt with in Lemma
  \ref{lem:special} below. 
\end{proof}
In case 0 of Lemma \ref{lem:cases}, the proof of the theorem is
trivial: $\nu_h(T_{m,n})= 0$, so $\rank\eta =\nu'_h(T_{m,n}) = 0$. We
assume from now on that $q \mid \frac{mn}{c}$, and use the method of
\S5 of \cite{CJKS}.  Let $\delta$ be the (oriented) $m$-string braid
$\sigma_{m-1}\sigma_{m-2}\cdots\sigma_1$. (The 3-string braids shown
in Figure \ref{fig:braids} should make our conventions sufficiently 
clear.) 
\setlength{\unitlength}{1pt}
\begin{figure}
\begin{center}
\begin{picture}(248,81)(10,45) {
\put(10,126){\line(0,-1){22}}
\put(10,104){\line(1,-1){10}}
\put(34,80){\line(-1,1){10}}
\put(34,80){\vector(0,-1){22}}
\put(34,126){\line(0,-1){22}}
\put(34,104){\line(-1,-1){24}}
\put(10,80){\vector(0,-1){22}}
\put(58,126){\vector(0,-1){68}}
\put(31,46){$\sigma_1$}
\put(110,126){\vector(0,-1){68}}
\put(134,126){\line(0,-1){22}}
\put(134,104){\line(1,-1){10}}
\put(158,80){\line(-1,1){10}}
\put(158,80){\vector(0,-1){22}}
\put(158,126){\line(0,-1){22}}
\put(158,104){\line(-1,-1){24}}
\put(134,80){\vector(0,-1){22}}
\put(131,46){$\sigma_2$}
\put(210,126){\line(0,-1){34}}
\put(210,92){\line(1,-1){10}}
\put(234,68){\line(-1,1){10}}
\put(234,68){\vector(0,-1){10}}
\put(234,126){\line(0,-1){10}}
\put(234,116){\line(1,-1){10}}
\put(258,92){\line(-1,1){10}}
\put(258,92){\vector(0,-1){34}}
\put(258,126){\line(0,-1){10}}
\put(258,116){\line(-1,-1){48}}
\put(210,68){\vector(0,-1){10}}
\put(226,46){$\sigma_2\sigma_1$}
}\end{picture}
\end{center}
\caption{}
\label{fig:braids}
\end{figure}

Let $\vec x = (x_1, \ldots, x_m) \in \Fld_p^m$, and represent
endomorphisms of $\Fld_p^m$ by $m\times m$ matrices acting on the
left. If the top arcs of $\delta$ are assigned the colors $x_1,
\ldots,x_m$ (from left to right), the colors on the bottom arcs are
determined by the vector $M_m(\vec x)$, where 
\[
  M_m=\left[\begin{matrix}
    0      & 0      & 0      &\cdots & 0      & 1 \\ 
    \theta & 0      & 0      &\cdots & 0      & 1 - \theta\\ 
    0      & \theta & 0      &\cdots & 0      & 1-\theta\\
    \vdots &\vdots  &\ddots  &\ddots &\vdots  & \vdots\\
    0      & 0      & \cdots &\theta & 0      & 1 -\theta\\
    0      & 0      & \cdots & 0     & \theta & 1-\theta
  \end{matrix}\right].
\]
(This is the transpose of the matrix in \cite{CJKS} because we use a
left action.)  Since $T_{m,n}$ is the closure of $\delta^n$, $\vec x$
determines a coloring of $T_{m,n}$ iff $\vec x$ is fixed by $M_m^n$,
so we may identify $\rep(T_{m,n},\Fld_p)$ with the space 
\[ 
  V_{m,n} =\{\,\vec x \in \Fld_p^m \mid M_m^n(\vec x) = \vec x\,\}.
\]
With this identification, the form $\eta$ is given as follows. Let
$\vec x$ and $\vec y$ be in $V_{m,n}$. Then 
\begin{equation}\label{eq:eta1}
  \eta(\vec x, \vec y) = \sum_{i = 0}^{n-1}\sum_{j=1}^{m-1}
  \cism\left(M_m^i(\vec x)_j\overline{M_m^i(\vec y)_m} - M_m^i(\vec x)_m
  \overline{M_m^i(\vec y)_j}\right).
\end{equation}
\begin{lemma}\label{lem:charpoly}\hfill
  \begin{enumerate}
    \item The characteristic polynomial $\det(T I_m - M_m)\in
      \Fld_p[T]$ of $M_m$ is equal to $(T-1)\frac{T^m - \theta^m}
      {T -\theta}$.   
    \item If $\vec x \in \Fld_p^m$ is an eigenvector of $M_m$, then
      $x_m \neq 0$. 
    \item The eigenspaces of $M_m$ are one-dimensional.
    \item If $\vec x \in \Fld_p^m$ is an eigenvector of $M_m$ of
    eigenvalue $\xi$, then $\vec x \in V_{m,n}$ iff $\xi^n = 1$. 
  \end{enumerate}
\end{lemma}
\begin{proof}
  Since each row of $M_m$ sums to 1, each row of $T I_m-M_m$ sums to
  $T-1$, and it follows that $\det(T I_m-M_m) = (T-1)\Delta_m$, where 
  \[
    \Delta_m=\det\left[\begin{matrix}
      T       & 0      & 0     & \cdots & 0     & 1 \\ 
      -\theta & T      & 0     & \cdots & 0     & 1\\
      0       &-\theta & T     &\cdots  & 0     & 1\\
      \vdots  &\vdots  &\ddots &\ddots  &\vdots &\vdots\\
      0       & 0      &\cdots &-\theta & T     & 1\\
      0       & 0      &\cdots & 0      & -\theta & 1
    \end{matrix}\right].
  \] 
  Expanding along the first row gives $\Delta_m = T\Delta_{m-1}
  +\theta^{m-1}$, and part (1) follows by induction. Suppose that
  $\vec x \in \Fld_p^m$ and $M_m(\vec x) = \xi \vec x$. Then 
  $\theta x_{i-1} + (1-\theta)x_m = \xi x_i$ for $2 \leq i \leq m$,
  and it follows that $x_m = 0$ implies $\vec x = 0$. This proves
  parts (2) and (3). Also $M_m^n(\vec x) = \xi^n\vec x$, which gives
  (4). 
\end{proof}
We make $\Fld_p^m$ into a module over $\Fld_p[T,T^{-1}]$ by setting
$T\vec x = M_m(\vec x)$; a generator for the order ideal is 
$\det(T I_m-M_m)$. Then $V_{m,n}$ is a submodule, and is annihilated
by the greatest common divisor of $\det(T I_m - M_m)$ and
$T^n-1$. Define $\maps{\sigma,\pi}{V_{m,n}}{\Fld_p}$ by 
$\sigma(\vec x) = \sum_{i=1}^m x_i$ and $\pi(\vec x) = x_m$. Note that
in the inner sum of \eqref{eq:eta1}, we may let $j$ range from 1 to
$m$ (as the extra terms cancel), so the equation may be rewritten as 
\begin{equation}\label{eq:eta2}
  \eta(\vec x, \vec y) = \cism\sum_{i = 0}^{n-1} 
  \left(\sigma(T^i\vec x)\overline{\pi(T^i\vec y)} 
  - \pi(T^i\vec x)\overline{\sigma(T^i\vec y)}\right).  
\end{equation}
Note also that
\begin{equation}\label{eq:sigmaT}
  \sigma(T\vec x) = \theta\sigma(\vec x) + m(1-\theta)x_m.
\end{equation}
For $f \in \Fld_p[T,T^{-1}]$, set $f^* = \bar f(T^{-1})$. Note that
$T$ preserves the form $\eta$, so that $\eta(f \vec x, \vec y) =
\eta(\vec x, f^*\vec y)$ for $f \in \Fld_p[T,T^{-1}]$ and $\vec x,
\vec y\in V_{m,n}$. 
\begin{lemma}\label{lem:orthog}
  Suppose that $\vec x, \vec y \in V_{m,n}$ have orders
  $f,g\in\Fld_p[T,T^{-1}]$, respectively, and that $f^*$ is coprime to
  $g$. Then $\eta(\vec x, \vec y) = 0$. 
\end{lemma}
\begin{proof} 
  This is a standard argument. Let $f_1, g_1\in\Fld_p[T,T^{-1}]$ be
  such that $f_1f^*+g_1g = 1$. Then 
  \[
    \eta(\vec x, \vec y)=\eta\bigl(\vec x,(f_1f^* + g_1g)\vec y\bigr)
    =\eta(f\vec x,f_1\vec y)+\eta(\vec x, g_1g\vec y)=0.
  \]
\end{proof}
We next prove Theorem \ref{thm:torus} in cases 1 and 2 of Lemma
\ref{lem:cases}. 
\begin{lemma}\label{lem:pdivm}
  If $q\mid \frac{mn}{c}$ and $p\mid m$, then $\rank\eta
  =\nu'_h(T_{m,n})$. 
\end{lemma}
\begin{proof} 
  Let $\vec x, \vec y \in V_{m,n}$. Since $p\mid m$, equation
  \eqref{eq:sigmaT} simplifies to $\sigma(T\vec x) = 
  \theta\sigma(\vec x)$, so \eqref{eq:eta2} gives  
  \begin{align*}
    \eta(\vec x, \vec y) &= \cism\sum_{i = 0}^{n-1} 
    \left(\theta^i\sigma(\vec x)\overline{\pi(T^i\vec y)} -
    \pi(T^i\vec x)\overline{\theta^i\sigma(\vec y)}\right)\\ 
    &= \cism\sum_{i = 0}^{n-1} 
    \left(\sigma(\vec x)\overline{\pi(\theta^{-i}T^i\vec y)} -  
    \pi(\theta^{-i}T^i\vec x)\overline{\sigma(\vec y)}\right)\\
    &= \cism\left(\sigma(\vec x)
    \overline{\pi\left(\theta^{-n+1}\frac{T^n-\theta^n}{ T-\theta}
    \vec y\right)} -
    \pi\left(\theta^{-n+1}\frac{T^n-\theta^n}{T-\theta}
    \vec x\right) \overline{\sigma(\vec y)}\right). 
  \end{align*}
  Suppose that $q \nmid n$, so that $\theta^n \neq 1$. Then 
  $\sigma(\vec x)=\sigma(T^n\vec x) = \theta^n\sigma(\vec x)$, and
  hence $\sigma(\vec x)=0$ for all $\vec x \in V_{m,n}$. In this case,
  $\rank\eta = 0$, as required. 

  Suppose then that $q \mid n$. Consider first case 1 of Lemma
  \ref{lem:cases}, where $p \mid c$ and $p\nmid \frac{m}{c}$. Then
  $\val_p(m) \leq \val_p(n)$. By Lemma \ref{lem:val}, $\val_{T-\theta}
  \bigl(\det(T I_m - M_m)\bigr) = p^{\val_p(m)}-1$, while $\val_{T-\theta}
  (T^n-1) = \val_{T-\theta}(T^n - \theta^n) = p^{\val_p(n)}$. It
  follows that $V_{m,n}$ is annihilated by
  $\frac{T^n-\theta^n}{T-\theta}$, and again we find $\rank \eta = 0$.  

  Finally, suppose $p \nmid n$. Let $V'$ be the submodule of $V_{m,n}$
  consisting of elements whose order is a power of $T-\theta$, and
  $V''$ the submodule of elements whose order is prime to
  $T-\theta$. Then $V_{m,n} = V' \oplus V''$, and by Lemma
  \ref{lem:orthog} this is an orthogonal sum (because $T-\theta$ and
  $(T-\theta)^*$ are equal up to units). Since $T-\theta$ divides
  $\det(T I_m-M_m)$, $M_m$ has an eigenvector of eigenvalue $\theta$
  in $\Fld_p^m$, and since $\theta^n=1$, this is in $V_{m,n}$. Thus
  $V'$ is non-zero. Because $T-\theta$ is not a repeated factor of
  $T^n-1$, every non-zero element of $V'$ is an eigenvector, so $V'$
  is one-dimensional. If $\vec x \in V''$ has order $f(T)$, then 
  $0 = \sigma\bigl(f(T)\vec x\bigr) = f(\theta)\sigma(\vec x)$, and  
  since $f(\theta)\neq 0$, $\sigma(\vec x)=0$. It follows that $\eta$
  is identically zero on $V''\times V''$, and hence that $\rank\eta$
  is 0 or 1 according as $\eta(\vec x, \vec x)$ is or is not 0 for
  $\vec x$ an eigenvector of eigenvalue $\theta$. For such $\vec x$,
  we have 
  \[
    \eta(\vec x, \vec x) = n\cism \left(\sigma(\vec x) \bar x_m - 
    x_m \overline{\sigma(\vec x)}\right). 
  \]
  For $1 \leq i < m$, $\theta x_{i+1} = (T\vec x)_{i+1} = \theta x_i
  +(1-\theta)x_m$, so $x_i - x_{i+1}= (1-\bar\theta)x_m$, and 
  $\sigma(\vec x) = \sum_{i=1}^{m-1} i(x_i - x_{i+1})
  =\frac{m(m-1)}{2}(1-\bar\theta)x_m$. If $p \neq 2$,
  $\frac{m(m-1)}{2}\equiv 0 \pmod p$, so $\sigma(\vec x) = 0$ and
  $\rank \eta = 0$. If $p=2$, we have $\bar\theta = \theta+1$, and so
  $\sigma(\vec x) = (m/2)\theta x_m$ and $\eta(\vec x, \vec x) =
  (m/2)n\cism x_m \bar x_m \in \fld_2$. Using part (2) of Lemma
  \ref{lem:charpoly}, this is non-zero iff $m \equiv 2\pmod 4$, and
  the proof is complete. 
\end{proof}
\begin{lemma}\label{lem:special}
  Suppose that $p=2$, $3 \mid n$, $m\equiv 2\pmod 4$, and 
  $2 \nmid n$. Then $\rad\eta = \rep'(T_{m,n}, \Fld_p)$.  
\end{lemma}
\begin{proof} 
  We have $\rank\eta = \nu'_h(T_{m,n})=1$, and it is enough to show
  that $\rad\eta \leq \rep'(T_{m,n})$. From the proof of Lemma
  \ref{lem:pdivm}, $\rad\eta$ is the kernel of $\sigma$.  From the
  proof of Lemma \ref{lem:repp}, it is enough to exhibit a single
  element $a$ of $A_2(T_{m,n})$ such that $ha=0$, $a \notin
  hA_2(T_{m,n})$, and $\vec x(a) = 0$ for all $\vec x\in \rad\eta$.
  In the diagram of $T_{m,n}$ corresponding to the braid $\delta^n$,
  let $a_{i,1},\ldots,a_{i,m}$ be the arcs at the top of the 
  $i^{\rm th}$ copy of $\delta$, for $ 1 \leq i \leq n$. Then
  $A_2(T_{m,n})$ is generated by $\{\,a_{i,j}\mid 1 \leq i \leq n, 1
  \leq j\leq m\}$, subject to the relations  
  \begin{alignat*}{2}
    a_{i,1} &= a_{i-1,m} && \quad\text{for $2\leq i \leq n$,}\\
    a_{i,j} &= t a_{i-1,j-1} + (1-t)a_{i-1,m} && \quad \text
                {for $2 \leq i \leq n$ and $2\leq j\leq m$,}\\   
    a_{1,1} &=a_{n,m}, && \\
    \text{and}\quad a_{1,j} &= t a_{n,j-1} + (1-t)a_{n,m} && \quad
                      \text{for $2\leq j\leq m$.}
  \end{alignat*}
  We may eliminate all the generators except $a_{1,1},\ldots,a_{1,m}$,
  and the relations between these may be expressed as follows.  Let
  $N_m$ be the $m\times m$ matrix over $\Lambda_2$ that results from
  replacing $\theta$ by $t$ in $M_m$, and $N_m^T$ its transpose. Then
  the relations are $a_{1,i}=\sum_{j=1}^m (N_m^n)_{ij}a_{1,j}$ for 
  $1 \leq i \leq m$. For $\vec f = (f_1,\ldots,f_m)\in \Lambda_2^m$,
  let $[\vec f\:] = \sum_{i=1}^mf_ia_{1,i}$. Then 
  $\bigl[(N_m^T)^n(\vec f\:)\bigr]= [\vec f\:]$.  Let $\vec f_0 =
  (1,1,\ldots,1)\in \Lambda_2^m$. Then (because $m$ is even)
  $N_m^T(\vec f_0) = t\vec f_0$, so $(t^n-1)[\vec f_0] = 0$. Since
  $3\mid n$ and $2 \nmid n$, we have $\val_h(t^n-1)=1$, so $t^n-1=gh$
  for some $g\in\Lambda_2$ which is coprime to $h$, and hence has
  $g(\theta)\neq 0$. Let $a = [g\vec f_0] \in A_2(T_{m,n})$. Then
  $ha=0$ and, for $\vec x \in V_{m,n}$, $\vec x(a) =
  g(\theta)(x_1+\cdots+x_m)$. Hence $\vec x(a)=0$ iff $\vec x \in \rad
  \eta$. Because there exists an $\vec x \in V_{m,n}$ with $\vec x(a)
  \neq 0$, $a \notin hA_2(T_{m,n})$, and we are done. 
\end{proof}
We let $\hat\Fld_p$ be the splitting field of $T^c-1$ over $\Fld_p$.
\begin{lemma}\label{lem:invol} 
  Suppose that $p \nmid c$ and $4\nmid \ord_c(p)$.  Then the
  involution $ x \mapsto \bar x$ of $\Fld_p$ extends uniquely to an
  involution of $\hat\Fld_p$. 
\end{lemma}
\begin{proof}
  Because $[\hat\Fld_p:\fld_p]$ is even, there is a unique non-trivial
  involution $\alpha$ of $\hat\Fld_p$.  Let $\mathbb K$ be the finite
  field of order $p^d$. Since $p\nmid c$, $T^c-1$ splits over 
  $\mathbb K$ iff $\mathbb K$ contains $c$ distinct $c^{\rm th}$ roots
  of unity. Since $\mathbb K^\times$ is cyclic of order $p^d-1$, this
  happens iff $c \mid p^d-1$; i.e., iff $\ord_c(p)\mid d$. Hence the
  splitting field $\hat\fld_p$ of $T^c-1$ over $\fld_p$ has
  $[\hat\fld_p:\fld_p] = \ord_c(p)$. If $\ord_c(p)$ is odd, then
  $[\hat\Fld_p:\Fld_p] = \ord_c(p)$, while if $\ord_c(p)$ is even then
  $\hat\fld_p = \hat\Fld_p$ and $[\hat\Fld_p:\Fld_p] = \ord_c(p)/2$. 
  Since $4 \nmid\ord_c(p)$, $[\hat\Fld_p:\Fld_p]$ is odd in either
  case. Hence the fixed field of $\alpha$ does not contain $\Fld_p$,
  so the restriction of $\alpha$ to $\Fld_p$ is non-trivial. Since 
  $x \mapsto \bar x$ is the only non-trivial automorphism of $\Fld_p$,
  the result follows. 
\end{proof}
From now on we assume that $q \mid \frac{mn}{c}$, $p \nmid mn$ and 
$4 \nmid \ord_c(p)$; in particular, Lemma \ref{lem:invol}
applies. That lemma having been established, no confusion will result
if we write the involution of $\hat\Fld_p$ as $x \mapsto \bar x$. We
consider the action of the matrix $M_m$ on $\hat\Fld_p^m$, and set  
\[
  \hat V_{m,n} = 
  \{\,\vec x \in \hat\Fld_p^m \mid M_m^n(\vec x) = \vec x\,\}.  
\]
Then $\hat V_{m,n}$ is the $\hat\Fld_p$-subspace of $\hat\Fld_p^m$
generated by $V_{m,n}$, and in particular has dimension
$\nu_h(T_{m,n})+1$. The action of $M_m$ makes $\hat \Fld_p^m$ into a
module over $\hat\Fld_p[T,T^{-1}]$ with $\hat V_{m,n}$ as a
submodule. Extend $\sigma$ and $\pi$ to 
$\maps{\sigma, \pi}{\hat V_{m,n}}{\hat\Fld_p}$ using the same  
formulas as before, and extend $\eta$ to a Hermitian form $\hat \eta$
on $\hat V_{m,n}$ given by the right-hand side of \eqref{eq:eta2}.   
Then $\rank\hat\eta =\rank\eta$, and Lemmas \ref{lem:charpoly} and
\ref{lem:orthog} hold with $\Fld_p$ replaced by $\hat\Fld_p$.  For 
$x \in \hat\Fld_p$, set $x^\bullet = \bar x^{-1}$; note that   
$x = x^\bullet$ iff $N(x)=1$.  
\begin{lemma} \label{lem:eigen} 
  Suppose that $q \mid \frac{mn}{c}$, $p \nmid mn$ and $4 \nmid
  \ord_c(p)$.  Then the number of distinct eigenvalues of $T | \hat
  V_{m,n}$ in $\hat\Fld_p$ is $\nu_h(T_{m,n})+1$, and if $\xi$ is an
  eigenvalue, so is $\xi^\bullet$. Further, $\theta$ is not an
  eigenvalue. 
\end{lemma}
\begin{proof} 
  Let $E$ be the set of eigenvalues; by Lemma \ref{lem:charpoly}, this
  is the set of common roots of the polynomials
  $(T-1)\frac{T^m-\theta^m}{T - \theta}$ and $T^n-1$ in
  $\hat\Fld_p$. Let $E'$ be the set of common roots of $T^m-\theta^m$
  and $T^n-1$. The symmetric difference of $E$ and $E'$ contains at
  most 1 and $\theta$. Since $p \nmid m$, $\theta$ is a root of
  $T^m-\theta^m$ but not of $\frac{T^m-\theta^m}{T- \theta}$; hence
  $\theta \notin E$, while $\theta \in E'$ iff $q \mid n$. On the
  other hand, $1 \in E$, while $1 \in E'$ iff $q \mid m$. Hence 
  \[
    |E|=\begin{cases}
      |E'|+1, & \text{if $q \nmid m$ and $q \nmid n$;}\\ 
      |E'|-1, & \text{if $q \mid m$ and $q \mid n$;}\\
      |E'|, &\text{otherwise.}
    \end{cases}
  \] 
  Let $c_1$ be the highest common factor of $c$ and $q$. Since 
  $q \mid \frac{mn}{c}$, we have $q = c_1m_1n_1$ where $c_1m_1\mid m$
  and $c_1n_1\mid n$, and then $m_1$ and $n_1$ are coprime. Let 
  $m = c_1m_1m_2$ and $n=c_1n_1n_2$.  By the Chinese Remainder
  Theorem, there is an integer $i$ with $i \equiv 0 \pmod {m_1}$ and
  $i \equiv 1 \pmod {n_1}$.  Set $\xi = \theta^i$. Since $ic_1n_1
  \equiv 0 \pmod{q}$, $\xi^n = (\theta^{ic_1n_1})^{n_2} = 1$, and
  since $ic_1m_1 \equiv c_1m_1 \pmod{q}$, $\xi^m =
  (\theta^{ic_1m_1})^{m_2} =(\theta^{c_1m_1})^{m_2} = \theta^m$. Thus
  $\xi \in E'$. Since $E'$ is non-empty, $|E'|$ is equal to the number
  of $c^{\rm th}$ roots of unity in $\hat\Fld_p$, which is $c$.  By
  case 3 of Lemma \ref{lem:cases}, we do have $|E| = \nu_h(T_{m,n}) +
  1$.  It was observed above that $\theta \notin E$. It is clear that
  $\xi\in E'$ iff $\xi^\bullet\in E'$; since $1^\bullet = 1$ and
  $\theta^\bullet=\theta$, the same is true of $E$, and we are done. 
\end{proof}
The next lemma gives all remaining cases of Theorem \ref{thm:torus}.
\begin{lemma} \label{lem:pndivmn} 
  Suppose that $q \mid \frac{mn}{c}$, $p \nmid mn$ and $4 \nmid
  \ord_c(p)$. Then $\rank\eta =\nu'_h(T_{m,n})$. 
\end{lemma}
\begin{proof} 
  By Lemma \ref{lem:eigen}, $\hat V_{m,n}$ is the direct sum of the
  eigenspaces of $T$, which are all one-dimensional, and the
  eigenvalues with norm different from 1 come in pairs $\xi$ and
  $\xi^\bullet$. An eigenvector of eigenvalue $\xi$ has order $T-\xi$,
  and $(T-\xi)^* = T^{-1}-\bar\xi = -T^{-1}\bar\xi(T-\xi^\bullet)$, so
  Lemma \ref{lem:orthog} shows that: 
  \begin{itemize}
    \item[(1)] an eigenspace corresponding to an eigenvalue of norm 1
      is orthogonal to the other eigenspaces; 
    \item[(2)] for a pair of eigenvalues $\xi$ and $\xi^\bullet$ of
      norm different from 1, the sum of the corresponding eigenspaces
      is orthogonal to all other eigenspaces, and $\hat\eta$ is zero
      on each of these spaces. 
  \end{itemize}
  Hence, $\rank\hat\eta$ is the number of eigenvalues $\xi$ such that
  $\hat\eta(\vec x, \vec y)\neq 0$ for eigenvectors $\vec x$ and 
  $\vec y$ with eigenvalues $\xi$ and $\xi^\bullet$. By case 3 of
  Lemma \ref{lem:cases}, $\nu'_h(T_{m,n})=\dim \hat V_{m,n}-1$, so we
  must show that there is just one eigenvalue for which this condition
  does not hold. 

  Let $\vec x$ be an eigenvector of eigenvalue $\xi$, and $\vec y$ an
  eigenvector of eigenvalue $\xi^\bullet$. From equation
  \eqref{eq:eta2},  
  \[
    \hat\eta(\vec x, \vec y) = n\cism\bigl(\sigma(\vec x)\bar y_m
    -x_m\overline{\sigma(\vec y)}\bigr).
  \]
  If $\xi = 1$, $\vec x$ and $\vec y$ are multiples of
  $(1,1,\ldots,1)$, so $\hat\eta(\vec x,\vec y)=0$.  Suppose $\xi \neq
  1$. Equation \eqref{eq:sigmaT} gives  
  \begin{align*}
    \xi\sigma(\vec x) &=\theta\sigma(\vec x) + m(1-\theta)x_m\\
    \text{and}\quad\xi^\bullet\sigma(\vec y) &= \theta\sigma(\vec y) +
      m(1-\theta)y_m.
  \end{align*}
  Hence
  \begin{align*}
    (\xi - \theta)\overline{(\xi^\bullet-\theta)}\hat\eta(\vec x, \vec y) 
    &=n\cism\bigl(\overline{(\xi^\bullet-\theta)}(\xi - \theta)
    \sigma(\vec x)\bar y_m - (\xi - \theta)x_m\overline {(\xi^\bullet
    - \theta)\sigma(\vec y)}\bigr)\\ 
    &=n\cism\bigl(\overline{(\xi^\bullet-\theta)}m(1 - \theta)x_m\bar y_m -
    (\xi - \theta)x_m\overline{m(1 - \theta)y_m}\bigr)\\ 
    &=mn\cism x_m\bar y_m\bigl(\bar\xi^\bullet\bar\theta(\theta-\xi)
    (1 - \theta) - (\xi -\theta)(1 - \bar\theta)\bigr)\\  
    &=mn\cism x_m\bar y_m(\xi - \theta)(1- \bar\theta)(\bar\xi^\bullet-1),
  \end{align*}
  which is non-zero by part (2) of Lemma \ref{lem:charpoly}. This
  completes the proof. 
\end{proof}
\section{Two-bridge knots}\label{sec:2bridge}
Let $M_+$ and $M_-$ be the $2\times 2$ matrices over $\Lambda$ given
by   
\[
  M_\pm=\left[\begin{matrix}
    t^{\pm1}&1-t^{\pm 1} \\ 
    t^{\pm 1}-1 & 2 -t^{\pm 1}
  \end{matrix}\right].
\]
\begin{lemma}\label{lem:mn}
  For $n\in \z$, we have
  \[
    M_\pm^n=\left[\begin{matrix}
      n(t^{\pm1}-1)+1&n(1-t^{\pm 1}) \\ 
      n(t^{\pm 1}-1) & n(1 -t^{\pm1})+1
    \end{matrix}\right].
  \]
\end{lemma}
\begin{proof} 
  The desired equation is equivalent to $M_\pm^n = nM_\pm-(n-1)I$. Since
  $M_\pm$ has determinant 1 and trace 2, $M_\pm^2-2M_\pm+I=0$. Using
  this in the form $M_\pm^2=2M_\pm-I$, the cases $n\geq 1$ follow by
  an upward induction, and a downward induction using
  $M_\pm^{-1}=2I-M_\pm$ gives the other cases. 
\end{proof}
Let $K(P,Q)$ be the two-bridge knot whose 2-fold branched cyclic cover
is the lens space $L(P,Q)$. Here $P$ is odd, $Q$ is coprime to $P$,
and we may assume that $0\leq Q< P$. Since the mirror image of
$K(P,Q)$ is $K(P,P-Q)$, we may further assume that $Q$ is even. Then 
$\frac{P}{Q}$ has a continued fraction expansion of the form
\begin{equation}\label{eq:cfrac}
  2m_1-\cfrac{1}{2n_1-\cfrac{1}{2m_2-\cdots -\cfrac{1}{2n_k}}}\,.
\end{equation}
Hence $K(P,Q)$ has the diagram shown in Figure \ref{fig:2bridge}, in
which each box represents the indicated number of full right-hand
twists. 
\setlength{\unitlength}{1pt}
\begin{figure}
\begin{center}
\begin{picture}(237,65)(40,25) {
\put(50,50){\framebox(30,30){}}
\put(60,63){$m_1$}
\put(100,70){\framebox(30,30){}}
\put(112,83){$n_1$}
\put(150,50){\framebox(30,30){}}
\put(160,63){$m_2$}
\put(217,70){\framebox(30,30){}}
\put(229,83){$n_k$}
\put(50,35){\line(1,0){140}}
\put(193,34.5){\dots}
\put(207,35){\line(1,0){40}}
\put(80,55){\line(1,0){70}}
\put(180,55){\line(1,0){10}}
\put(193,54.5){\dots}
\put(207,55){\line(1,0){40}}
\put(80,75){\line(1,0){20}}
\put(130,75){\line(1,0){20}}
\put(180,75){\line(1,0){10}}
\put(193,74.5){\dots}
\put(207,75){\line(1,0){10}}
\put(50,95){\line(1,0){50}}
\put(130,95){\line(1,0){60}}
\put(193,94.5){\dots}
\put(207,95){\line(1,0){10}}
\put(50,45){\oval(20,20)[l]}
\put(50,85){\oval(20,20)[l]}
\put(247,65){\oval(20,20)[r]}
\put(247,65){\oval(60,60)[r]}
}\end{picture}
\end{center}
\caption{}
\label{fig:2bridge}
\end{figure}

If $D$ is any oriented tangle diagram, we may define $A(D)$ and
$A_p(D)$ to be the modules over $\Lambda$ and $\Lambda_p$,
respectively,  with generators the arcs of $D$ and a relation for each
crossing just as for link diagrams. Then, setting
$\rep(D,\Fld_p)=\Hom_{\Lambda_p}(A_p(D),\Fld_p)$, we may define
$B(f)\in \fld_p$ for $f\in \rep(D,\Fld_p)$ as before. 
\begin{lemma}\label{lem:1box1}
  Let $D$ be the tangle diagram of Figure \ref{fig:1box}(a). In
  $A(D)$, we have the relations 
  \begin{align*}
    c &= n(t-1)(a-b) + a\\
    \text{and}\quad d &= n(t-1)(a-b) +b.
  \end{align*}
  For $f\in\rep(D,\Fld_p)$, $B(f) = n\cism(\bar\theta-\theta)
  N\bigl(f(a-b)\bigr)$.  
\end{lemma}
\setlength{\unitlength}{1pt}
\begin{figure}
\begin{center}
\begin{picture}(164,52)(33,28) {
\put(50,50){\framebox(30,30){}}
\put(62,63){$n$}
\put(50,55){\vector(-1,0){10}}
\put(40,75){\vector(1,0){10}}
\put(90,55){\vector(-1,0){10}}
\put(80,75){\vector(1,0){10}}
\put(33,52){$a$}
\put(33,72){$b$}
\put(93,52){$c$}
\put(93,72){$d$}
\put(59,30){(a)}
\put(150,50){\framebox(30,30){}}
\put(162,63){$n$}
\put(140,55){\vector(1,0){10}}
\put(150,75){\vector(-1,0){10}}
\put(180,55){\vector(1,0){10}}
\put(190,75){\vector(-1,0){10}}
\put(133,52){$a$}
\put(133,72){$b$}
\put(193,52){$c$}
\put(193,72){$d$}
\put(159,30){(b)}
}\end{picture}
\end{center}
\caption{}
\label{fig:1box}
\end{figure}

\begin{proof}
  By Lemma \ref{lem:mn}, the claimed relations in $A(D)$ are
  equivalent to the matrix equation 
  \[
    \left[\begin{matrix}
      c\\ d
    \end{matrix}\right]
    =M_+^n\left[\begin{matrix}
      a\\ b
    \end{matrix}\right].
  \]
  It is enough to check the cases $n=\pm 1$, and this is left to the
  reader. For the second part, let the images of $a$, $b$, $c$ and $d$
  under $f$ be $\alpha$, $\beta$, $\gamma$ and $\delta$. When $n=1$,
  both crossings are negative, and we have
  $B(f)=-\phi(\alpha,\beta)-\phi(\delta,\gamma)
  =\phi(\gamma,\delta)-\phi(\alpha,\beta)$. When $n=-1$, the crossings
  are positive, and we have
  $B(f)=\phi(\beta,\alpha)+\phi(\gamma,\delta)
  =\phi(\gamma,\delta)-\phi(\alpha,\beta)$. Hence
  $B(f)=\phi(\gamma,\delta)-\phi(\alpha,\beta)$ for any $n$. From the
  first part, 
  \begin{align*}
    \gamma\bar\delta - \bar\gamma\delta
    &=n(\theta-1)(\alpha-\beta)(\bar\beta-\bar\alpha)
    +n(\bar\theta-1)(\bar\alpha-\bar\beta)(\alpha-\beta)
    +\alpha\bar\beta-\bar\alpha\beta  \\
    &=n(\bar\theta-\theta)N(\alpha-\beta)
    +\alpha\bar\beta-\bar\alpha\beta.
  \end{align*}
  Hence
  \begin{align*}
    \phi(\gamma,\delta)-\phi(\alpha,\beta)
    &=\cism(\gamma\bar\delta - \bar\gamma\delta
    -\alpha\bar\beta+\bar\alpha\beta)\\
    &=n\cism(\bar\theta-\theta)N(\alpha-\beta),
  \end{align*}
  as required.
\end{proof}
\begin{lemma}\label{lem:1box2}
  Let $D$ be the tangle diagram of Figure \ref{fig:1box}(b). In
  $A(D)$, we have the relations 
  \begin{align*}
    c &= n(t^{-1}-1)(a-b) + a\\
    \text{and}\quad d &= n(t^{-1}-1)(a-b) +b.
  \end{align*}
  For $f\in\rep(D,\Fld_p)$, $B(f) =0$.
\end{lemma}
\begin{proof}
 The first part is proved as in Lemma \ref{lem:1box1}, using $M_-$ in
  place of $M_+$. In the Boltzmann weight, the contributions from the
  crossings cancel in pairs. 
\end{proof}
We set $\nabla = -(t-1)(t^{-1}-1) = t - 2 + t^{-1}\in \Lambda$.
\begin{lemma}\label{lem:2box}
  Let $D$ be the tangle diagram of Figure \ref{fig:2box}. In $A(D)$,
  we have the relations 
  \[
    \left[\begin{matrix}a'-b' \\ b'-c'\end{matrix}\right]
    =\left[\begin{matrix}
      1+mn\nabla&-n(t^{-1}-1)\\m(t-1)&1\end{matrix}\right]
    \left[\begin{matrix}a-b \\ b-c\end{matrix}\right].
  \]
\end{lemma}
\setlength{\unitlength}{1pt}
\begin{figure}
\begin{center}
\begin{picture}(116,50)(33,50) {
\put(50,50){\framebox(30,30){}}
\put(60,63){$m$}
\put(100,70){\framebox(30,30){}}
\put(112,83){$n$}
\put(33,52){$a$}
\put(50,55){\vector(-1,0){10}}
\put(140,55){\vector(-1,0){60}}
\put(143,52){$a'$}
\put(33,72){$b$}
\put(40,75){\vector(1,0){10}}
\put(80,75){\vector(1,0){20}}
\put(87,65){$d$}
\put(130,75){\vector(1,0){10}}
\put(143,72){$b'$}
\put(33,92){$c$}
\put(100,95){\vector(-1,0){60}}
\put(140,95){\vector(-1,0){10}}
\put(143,92){$c'$}
}\end{picture}
\end{center}
\caption{}
\label{fig:2box}
\end{figure}

\begin{proof}
  From Lemmas \ref{lem:1box1} and \ref{lem:1box2} we have
  \begin{align*}
    a'&= m(t-1)(a-b)+a,\\
    d\phantom{'} &= m(t-1)(a-b)+b,\\
    b'&= n(t^{-1}-1)(d-c)+d,\\
    \text{and}\quad c'&= n(t^{-1}-1)(d-c)+c.
  \end{align*}
  The result follows.
\end{proof}
Given integers $m_1, n_1,\ldots,m_k,n_k$, we define polynomials
$\alpha_i$, $\beta_i$, $\gamma_i$ and $\delta_i$ in $\z[x]$ by
$\alpha_0=\delta_0 = 1$, $\beta_0=\gamma_0=0$, and, 
for $1\leq i\leq k$,  
\begin{align*}
  \alpha_i &= (1+m_in_ix)\alpha_{i-1} + n_ix\gamma_{i-1},\\
  \beta_i &= (1+m_in_ix)\beta_{i-1}+n_i\delta_{i-1},\\
  \gamma_i &= m_i\alpha_{i-1}+\gamma_{i-1},\\
  \text{and}\quad\delta_i &=m_i x\beta_{i-1}+\delta_{i-1}.
\end{align*}
\begin{lemma}\label{lem:2kbox}
  Let $D$ be the tangle diagram of Figure \ref{fig:2kbox}. In $A(D)$,
  we have the relations 
  \begin{equation}\label{eq:2kbox}
    \left[\begin{matrix}a_i-b_i \\ b_i-c_i\end{matrix}\right]
    =\left[\begin{matrix}
      \alpha_i(\nabla) & -(t^{-1}-1)\beta_i(\nabla)\\
      (t-1)\gamma_i(\nabla)&\delta_i(\nabla)
    \end{matrix}\right]
    \left[\begin{matrix}a_0-b_0 \\ b_0-c_0\end{matrix}\right]
  \end{equation}
  for $0\leq i\leq k$.
\end{lemma}
\setlength{\unitlength}{1pt}
\begin{figure}
\begin{center}
\begin{picture}(240,56)(29,44) {
\put(50,50){\framebox(30,30){}}
\put(60,63){$m_1$}
\put(100,70){\framebox(30,30){}}
\put(112,83){$n_1$}
\put(150,50){\framebox(30,30){}}
\put(160,63){$m_2$}
\put(217,70){\framebox(30,30){}}
\put(229,83){$n_k$}
\put(29,52){$a_0$}
\put(50,55){\vector(-1,0){10}}
\put(150,55){\vector(-1,0){70}}
\put(137,45){$a_1$}
\put(190,55){\vector(-1,0){10}}
\put(193,54.5){\dots}
\put(257,55){\vector(-1,0){50}}
\put(260,52){$a_k$}
\put(29,72){$b_0$}
\put(40,75){\vector(1,0){10}}
\put(80,75){\vector(1,0){20}}
\put(130,75){\vector(1,0){20}}
\put(137,65){$b_1$}
\put(180,75){\vector(1,0){10}}
\put(193,74.5){\dots}
\put(207,75){\vector(1,0){10}}
\put(247,75){\vector(1,0){10}}
\put(260,72){$b_k$}
\put(29,92){$c_0$}
\put(100,95){\vector(-1,0){60}}
\put(190,95){\vector(-1,0){60}}
\put(137,85){$c_1$}
\put(193,94.5){\dots}
\put(217,95){\vector(-1,0){10}}
\put(257,95){\vector(-1,0){10}}
\put(260,92){$c_k$}
}\end{picture}
\end{center}
\caption{}
\label{fig:2kbox}
\end{figure}

\begin{proof} 
  Using Lemma \ref{lem:2box} and the fact that the matrix  
  \[
    \left[\begin{matrix}
      \alpha_i(\nabla) & -(t^{-1}-1)\beta_i(\nabla)\\
      (t-1)\gamma_i(\nabla)&\delta_i(\nabla)
    \end{matrix}\right]
  \]
  is equal to
  \[
    \left[\begin{matrix} 
      1+m_in_i\nabla&-n_i(t^{-1}-1)\\ 
      m_i(t-1)&1
    \end{matrix}\right]
    \left[\begin{matrix}
      \alpha_{i-1}(\nabla) & -(t^{-1}-1)\beta_{i-1}(\nabla)\\
      (t-1)\gamma_{i-1}(\nabla)&\delta_{i-1}(\nabla)
    \end{matrix}\right],
  \]
  this follows by induction.
\end{proof}
\begin{lemma}\label{lem:2bridge}
  Let $K$ be the two-bridge knot $K(P,Q)$, where $\frac{P}{Q}$ has the
  continued fraction expansion \eqref{eq:cfrac}. Then $A(K)$ has a
  presentation with generators $a$ and $d$ and one relation
  $\alpha_k(\nabla)d=0$. For $f\in \rep(K,\Fld_p)$,  
  \[
    B(f) = \cism(\bar\theta-\theta)N\bigl(f(d)\bigr)
    \sum_{i=1}^km_i\alpha_{i-1}\bigl(\nabla(\theta)\bigr)^2.
  \] 
\end{lemma} 
In the expression $\alpha_{i-1}\bigl(\nabla(\theta)\bigr)$, and
similar expressions occuring later, it is to be understood that the
coefficients of the polynomials are taken modulo $p$. 
\begin{proof} 
  Let $D$ be the tangle diagram of Figure \ref{fig:2kbox}. It is clear
  that $A(D)$ is the free $\Lambda$-module on $a_0$, $b_0$ and
  $c_0$. Let $c_k = f_1a_0+f_2b_0+f_3c_0$, where $f_1,f_2,f_3\in
  \Lambda$. Then similar expressions for $a_k$ and $b_k$ are
  determined by the case $i=k$ of the relations \eqref{eq:2kbox}. Put
  differently, the initial presentation (given by all arcs and
  crossings) of $A(D)$ is equivalent to a presentation with generators
  $a_0$, $b_0$, $c_0$, $a_k$, $b_k$ and $c_k$, and the relations 
  \begin{equation}\label{eq:rel1}
    \begin{split}
      a_k-b_k &=
      \alpha_k(\nabla)(a_0-b_0)-(t^{-1}-1)\beta_k(\nabla)(b_0-c_0),\\
      b_k-c_k
      &=(t-1)\gamma_k(\nabla)(a_0-b_0)+\delta_k(\nabla)(b_0-c_0),\\
      \text{and}\quad\phantom{b_k-{}} c_k &=f_1a_0+f_2b_0+f_3c_0. 
    \end{split}
  \end{equation}
  The diagram of $K$ in Figure \ref{fig:2bridge} is obtained from $D$
  by joining the free ends without introducing any more
  crossings. Consider instead the diagram where one curl is introduced
  in each connecting arc. This gives a presentation of $A(K)$ which is
  obtained from the initial presentation of $A(D)$ by adding a
  relation for each of the three extra crossings. However, in this
  presentation, any relation is a consequence of the others, so we
  need only add 
  \begin{equation}
    b_0 =c_0,\quad a_k=b_k.\label{eq:rel2}
  \end{equation}
  Hence $A(K)$ has a a presentation with generators $a_0$, $b_0$,
  $c_0$, $a_k$, $b_k$ and $c_k$, and the relations \eqref{eq:rel1} and
  \eqref{eq:rel2}. In the presence of relations \eqref{eq:rel2},
  relations \eqref{eq:rel1} are equivalent to 
  \begin{align*}
    0 &= \alpha_k(\nabla)(a_0-b_0),\\
    a_k-c_k &=(t-1)\gamma_k(\nabla)(a_0-b_0),\\
    \text{and}\quad\phantom{a_k-{}} c_k &=f_1a_0+(f_2+f_3)b_0. 
  \end{align*}
  We may now eliminate four of the generators to arrive at a
  presentation of $A(K)$ with generators $a_0$ and $b_0$ and the
  single relation $\alpha_k(\nabla)(a_0-b_0)=0$. Setting $a=a_0$ and
  $d=a_0-b_0$ gives the claimed presentation. 

  By Lemmas (\ref{lem:1box1}) and (\ref{lem:1box2}), $B(f)=
  \sum_{i=1}^k m_i\cism(\bar\theta-\theta)
  N\bigl(f(a_{i-1}-b_{i-1})\bigr)$. In $A(K)$, we have
  $a_i-b_i=\alpha_i(\nabla)d$, so
  $f(a_i-b_i)=\alpha_i\bigl(\nabla(\theta)\bigr)f(d)$ and (since
  $\nabla(\theta)\in \fld_p$) $N\bigl(f(a_i-b_i)\bigr)
  =\alpha_i\bigl(\nabla(\theta)\bigr)^2N\bigl(f(d)\bigr)$, giving the
  required formula for $B(f)$.  
\end{proof}
\begin{theorem}\label{thm:2bridge}
  Let $K$ be a two-bridge knot. Then $\rad\eta=\rep'(K,\Fld_p)$.
\end{theorem} 
\begin{proof} 
  Recall that $h=t^2+\kappa t + 1$ for some $\kappa\in\fld_p$. Set
  $\lambda=\nabla(\theta)= -\kappa-2$, and note that $\lambda \neq 0$
  because $h$ is irreducible. Suppose that $K$ is determined by the
  continued fraction \eqref{eq:cfrac}. Since $\nu_h(K) \leq 1$ and
  $\nu_0(K)=1$, Lemma \ref{lem:c1impc2} shows that we need only prove
  that $\rank\eta=\nu_h'(K)$. Since $\rank\eta\leq \nu_h(K)$ and
  $\nu_h'(K)\leq \nu_h(K)$, this is trivially true if
  $\nu_h(K)=0$. Suppose for the rest of the proof that $\nu_h(K)=1$;
  by Lemma \ref{lem:2bridge}, this is equivalent to
  $\alpha_k(\lambda)=0$. It is enough to show that $\nu_h'(K)=0$ iff
  $\rank\eta=0$. We have $\nu_h'(K)=0$ iff $\lambda$ is a repeated
  root of $\alpha_k$, which is true iff $\lambda$ is a root of the
  derivative $\dot\alpha_k$. Also, $\rank\eta = 0$ iff $B(f)=0$ for
  all $f\in\rep(K,\Fld_p)$, and by Lemma \ref{lem:2bridge} again, this
  is true iff $\epsilon_k(\lambda)=0$ where
  $\epsilon_k=\sum_{i=1}^km_i\alpha_{i-1}^2\in\z[x]$. Thus we are
  reduced to proving that $\dot\alpha_k(\lambda) = 0$ iff
  $\epsilon_k(\lambda)=0$.  

  If $p\mid n_k$, we may replace $n_k$ by 0 without changing
  $\alpha_k$ or $\epsilon_k$ modulo $p$, and then the new $K$ has a
  shorter continued fraction. Thus we may assume that $p\nmid n_k$. We
  have $\alpha_k\delta_k - x\beta_k\gamma_k=1$, as one may see
  directly from the definitions, or more easily by using the matrix
  form of the recursion from the proof of Lemma \ref{lem:2kbox}. Hence
  $\gamma_k(\lambda)\neq 0$. From the definitions, we have $\alpha_k =
  \alpha_{k-1} +n_k x \gamma_k$, and hence $\alpha_{k-1}(\lambda)\neq 
  0$. We also have $\alpha_k=\alpha_{k-1}+n_k x\sum_{i=1}^k
  m_i\alpha_{i-1}$. Using this recursive relation, we prove in an
  appendix that $\dot\alpha_k\alpha_{k-1}-\alpha_k\dot\alpha_{k-1} =
  n_k \epsilon_k$, and it follows that $\dot\alpha_k(\lambda)=0$ iff
  $\epsilon_k(\lambda)=0$, completing the proof. 
\end{proof}
\appendix
\section{A polynomial identity}
It will be convenient to replace the integers $m_1,n_1,\ldots$ of
\S\ref{sec:2bridge} by indeterminates. Thus, let $R$ be the ring of
integer polynomials in variables $x,y_1,z_1,y_2,z_2,\ldots\,$. Define
polynomials $f_0, f_1, f_2,\ldots\in R$ inductively by $f_0=1$ and   
\[
  f_k = f_{k-1}+xz_k\sum_{i=1}^k y_i f_{i-1}
\] 
for $k \geq 1$. Comparing this with the recursion for the polynomials
$\alpha_i$ from the proof of Theorem \ref{thm:2bridge}, one sees that
Lemma \ref{lem:polyiden} below is what is needed to complete that
proof. 

Let $M\subset R$ be the set of monomials of the form  
\[
  \mu =y_{i_1}z_{j_1}y_{i_2}z_{j_2}\ldots y_{i_l}z_{j_l},\quad
  1\leq i_1 \leq j_1 < i_2 \leq j_2 < \cdots < i_l \leq j_l.
\]
We define the length of $\mu$ to be $l(\mu)=l$, and its height to be 0
if $l=0$ (when $\mu=1$) and $j_l$ if $l > 0$. For $k \geq 1$, let
$M_k$ be the set of elements of $M$ of height less than $k$. 
\begin{lemma}\label{lem:monom}
  \begin{align}
    \text{For $k\geq 0$,}\quad f_k &=\sum_{\mu\in M_{k+1}}x^{l(\mu)}\mu;
              \label{eq:fk1}\\
    \text{for $k \geq 1$,}\quad f_k &=\sum_{\mu\in M_{k}}x^{l(\mu)}\mu
    +\sum_{i=1}^k\sum_{\mu\in M_{i}}x^{l(\mu)+1}\mu y_i z_k.\label{eq:fk2}
  \end{align}
\end{lemma}
\begin{proof}
  Clearly \eqref{eq:fk1} holds for $k=0$. For $k\geq 1$, the
  right-hand sides of \eqref{eq:fk1} and \eqref{eq:fk2} are equal
  because $M_{k+1}$ is the disjoint union of $M_k$ and the sets
  $\{\,\mu y_i z_k\mid \mu \in M_i\,\}$ for $1\leq i\leq k$. The
  result follows by induction on $k$. 
\end{proof}
\begin{lemma}\label{lem:polyiden}
  Let $\dot f_k$ be the derivative of $f_k$ with respect to $x$. Then,
  for $k\geq 1$, 
  \[
    \dot f_kf_{k-1}-f_k \dot f_{k-1} = z_k\sum_{i=1}^k y_i f_{i-1}^2.
  \]
\end{lemma}
\begin{proof}
  Let $S=\{\,(i,\mu,\nu)\mid 1 \leq i < k, \mu \in M_i, \nu\in
  M_k-M_i\,\}$. By Lemma \ref{lem:monom}, we have  
  \begin{align*}
    \dot f_k f_{k-1}
    &= \Bigl(\sum_{\mu\in M_{k}}l(\mu)x^{l(\mu)-1}\mu
    +\sum_{i=1}^k\sum_{\mu\in M_{i}}\bigl(l(\mu)+1\bigr)x^{l(\mu)}\mu y_i
    z_k \Bigr)\Bigl(\sum_{\nu\in M_k}x^{l(\nu)}\nu\Bigr)\\ 
    &=\sum_{\mu,\nu\in M_{k}}l(\mu)x^{l(\mu)+l(\nu)-1}\mu\nu +\sum_{i=1}^k
    \sum_{\substack{\mu\in M_{i} \\ \nu\in M_k}}
    \bigl(l(\mu)+1\bigr)x^{l(\mu)+l(\nu)}\mu\nu y_i z_k\\ 
    &=\sum_{\mu, \nu\in M_{k}}l(\mu)x^{l(\mu)+l(\nu)-1}\mu\nu
    +\sum_{i=1}^k\sum_{\mu,\nu\in M_i}
    \bigl(l(\mu)+1\bigr)x^{l(\mu)+l(\nu)}\mu\nu y_i z_k\\
    &\qquad{}+\sum_{(i,\mu,\nu)\in S}
    \bigl(l(\mu)+1\bigr)x^{l(\mu)+l(\nu)}\mu\nu y_i z_k.
  \end{align*}
  Similarly, 
  \begin{align*}
    f_k\dot f_{k-1} &= \Bigl(\sum_{\mu'\in M_{k}}
    x^{l(\mu')}\mu' +\sum_{i=1}^k\sum_{\mu'\in M_{i}}
    x^{l(\mu')+1}\mu' y_i z_k \Bigr)\Bigl(\sum_{\nu'\in M_k}
    l(\nu')x^{l(\nu')-1}\nu'\Bigr)\\ 
    &=\sum_{\mu, \nu\in M_{k}}l(\mu)x^{l(\mu)+l(\nu)-1}\mu\nu 
    +\sum_{i=1}^k\sum_{\mu,\nu\in M_i} l(\mu)x^{l(\mu)+l(\nu)}\mu\nu y_i z_k\\
    &\qquad{}+\sum_{(i',\mu',\nu')\in S}
    l(\nu')x^{l(\mu')+l(\nu')}\mu'\nu' y_{i'} z_k. 
  \end{align*}
  (In the first two sums of the final expression, we have set
  $\mu=\nu'$ and $\nu=\mu')$. Let $(i,\mu,\nu)\in S$. Then $\nu$ can
  be written uniquely in the form $\nu=\mu' y_{i'}z_j$ with $1\leq
  i'\leq j<k$ and $\mu' \in M_{i'}$. Since $\nu \notin M_i$, we also
  have $i\leq j$. Set $\nu' = \mu y_i z_j$. Then $\nu'\in M_k-M_{i'}$,
  so $(i',\mu', \nu')\in S$. The function $(i,\mu,\nu)\mapsto
  (i',\mu', \nu')$ is an involution of $S$. Since $l(\mu') =
  l(\nu)-1$, $l(\nu')=l(\mu)+1$, and $\mu'\nu' y_{i'}=\mu\nu y_i$, it
  follows that 
  \[
    \sum_{(i,\mu,\nu)\in S} \bigl(l(\mu)+1\bigr)x^{l(\mu)+l(\nu)}
    \mu\nu y_i z_k =\sum_{(i',\mu',\nu')\in S}
    l(\nu')x^{l(\mu')+l(\nu')}\mu'\nu' y_{i'} z_k.
  \] 
  Hence 
  \begin{align*} 
    \dot f_k f_{k-1}-f_k \dot f_{k-1} &=
    \sum_{i=1}^k\sum_{\mu,\nu\in M_i}x^{l(\mu)+l(\nu)}\mu\nu y_i z_k\\
    &=z_k\sum_{i=1}^k\Bigl(\sum_{\mu\in M_i} x^{l(\mu)}\mu\Bigr)^2 y_i\\
    &=z_k\sum_{i=1}^k y_i f_{i-1}^2,
  \end{align*}
  as required.
\end{proof}

\end{document}